\documentclass[11pt]{article}
\usepackage{amsmath,amssymb,amsfonts}
\usepackage{geometry}
\usepackage{graphicx}
\usepackage{bm}
\title{Melnikov Analysis of Deterministic and Stochastic Manifold Splitting in the Kuramoto--Sivashinsky Equation}
\author{Sumita Datta$^{1,2}$\\
{\small $^{1}$Department of Pure and Applied Mathematics, Alliance University,}\\
{\small Bengaluru 562 106, India}\\
{\small $^{2}$Department of Physics, University of Texas at Arlington,}\\
{\small Texas 76019, USA}}
\date{\today}

\begin{document}
\maketitle

\begin{abstract}
We develop a Melnikov framework for the Kuramoto--Sivashinsky (KS) equation under weak deterministic and stochastic forcing. By treating KS as an infinite-dimensional dynamical system, we derive a Melnikov functional that measures splitting of stable and unstable manifolds of a homoclinic orbit. Periodic forcing leads to phase-dependent transverse intersections, while stochastic forcing produces random manifold splitting characterized by a variance determined by the adjoint solution. This provides a geometric mechanism linking invariant manifold theory to spatiotemporal chaos in dissipative partial differential equations.
\end{abstract}

\newpage
\section{Introduction}

The Kuramoto--Sivashinsky (KS) equation
\[
u_t + u_{xx} + u_{xxxx} + u u_x = 0
\]
is a canonical model for pattern formation and spatiotemporal chaos in dissipative extended systems. It arises in diverse physical contexts including flame front propagation, interfacial instabilities in thin films, reaction–diffusion systems, and plasma turbulence \cite{Kuramoto1984,Sivashinsky1977,CrossHohenberg1993, EckmannWayne1991}. Despite its deceptively simple form, the KS equation exhibits a rich hierarchy of behaviors ranging from steady patterns to low-dimensional chaos and fully developed spatiotemporal turbulence.

A key feature of the KS dynamics is the presence of coherent structures embedded in chaotic states. Extensive numerical and theoretical work has shown that the attractor of the KS flow possesses low-dimensional features governed by inertial manifold structure and unstable coherent solutions \cite{Temam1997, Robinson2001,FoiasSellTemam1988,ImkellerSchmalfuss2001}. Computer-assisted proofs have rigorously established chaotic dynamics in KS through finite-dimensional reductions and symbolic dynamics constructions \cite{Wilczak2003}. These studies suggest that, despite the infinite-dimensional phase space, essential dynamical mechanisms may be captured by geometric structures such as invariant manifolds and their intersections.

In finite-dimensional systems, one of the classical routes to chaos is homoclinic manifold splitting detected by Melnikov theory \cite{Wiggins2003,GuckenheimerHolmes1983,Smale1965}. Transverse intersections between stable and unstable manifolds of a saddle equilibrium generate homoclinic tangles and Smale horseshoes, providing a geometric explanation for chaotic dynamics. Extending this theory to dissipative partial differential equations has been an important challenge in infinite-dimensional dynamical systems \cite{Henry1981, BatesLuZeng1998}. For PDEs like KS, the phase space is a function space, invariant manifolds become infinite-dimensional objects, and classical ODE tools must be reformulated in a functional-analytic setting.

Previous studies of KS chaos have largely focused on numerical observation, inertial manifolds, and unstable periodic orbits \cite{Wilczak2003}. However, a direct Melnikov-type geometric mechanism for manifold splitting in KS, particularly under deterministic and stochastic perturbations\cite{ArnoldCrauelWihstutz1999,Schmalfuss1997,LuSchmalfuss2000}, has not been systematically developed. Understanding how weak forcing or noise perturbs homoclinic structures is essential for explaining transitions from coherent dynamics to turbulence \cite{HolmesLumleyBerkooz1996,Cvitanovic2005,Viswanath2007, Aubry1988} in extended systems.

In this work we develop a Melnikov framework \cite{ChowHale1982,Wiggins1990,Gelfreich2000}  for the Kuramoto–Sivashinsky equation treated as an infinite-dimensional dynamical system. We derive a Melnikov functional based on the adjoint linearized equation along a homoclinic orbit and show how it measures splitting of stable and unstable manifolds under weak perturbations. Both periodic and stochastic forcing are considered. Periodic forcing leads to phase-dependent transverse intersections, while stochastic forcing produces random manifold separation whose variance is determined by the adjoint solution.

The main contributions of this paper are:

\begin{itemize}
\item Extension of classical Melnikov analysis to a dissipative PDE setting;
\item Derivation of an adjoint-based Melnikov functional for the KS equation;
\item A geometric interpretation of deterministic and noise-driven spatiotemporal chaos;
\item A bridge between invariant manifold theory, stochastic PDEs, and numerical KS turbulence.
\end{itemize}

This work provides a unified geometric mechanism explaining how weak perturbations trigger complex dynamics in dissipative extended systems.
The Kuramoto--Sivashinsky (KS) equation is a prototypical dissipative partial differential equation exhibiting spatiotemporal chaos and low-dimensional coherent structures embedded within turbulent dynamics \cite{CrossHohenberg1993}.

A central question is identifying dynamical mechanisms responsible for the transition from coherent behavior to chaos. In finite-dimensional systems, homoclinic manifold splitting detected via Melnikov theory provides a classical route to chaos \cite{Wiggins1990,GuckenheimerHolmes1983}. Extending such ideas to dissipative infinite-dimensional systems remains an active area of research \cite{Temam1997, Robinson2001}.
In this work, we derive a Melnikov functional for KS under weak deterministic or stochastic forcing. This connects invariant manifold splitting to the onset of spatiotemporal chaos.

In Section 1, we introduce the problem. Section 2 provides a dynamical systems concepts in Melnikov analysis. Section 3 and 4 contain the derivation of Melnikov functional for KS equation subject to periodic and and stochastic forcing respectively. In Section 5 and 6 we have the formulas and parameters used in the numerical calculations of Melnokov functions. In section 7 we present our results. We conclude the paper in Section 8. Appendices A and B describe the classical Melnokov theory and derivation of adjoint equations.  

\section{Background on Invariant Manifolds and Homoclinic Orbits}

This section provides the dynamical systems concepts underlying the Melnikov analysis.

\subsection{Dynamical systems viewpoint}

A partial differential equation such as the Kuramoto--Sivashinsky (KS) equation can be viewed as an ordinary differential equation in an infinite-dimensional phase space. If we denote

\[
u_t = \mathcal{K}(u),
\]

then $u(t)$ evolves in a function space $H$ (e.g., $L^2_{\text{per}}$). Each function $u(x)$ represents a point in phase space, and the KS flow generates a trajectory $u(t)$ in this space.

Thus the KS equation is a dynamical system:

\[
\frac{d}{dt} u(t) = \mathcal{K}(u(t)).
\]

\subsection{Equilibria and their stability}

A steady solution $u_s(x)$ satisfies

\[
\mathcal{K}(u_s) = 0.
\]

Such a state corresponds to a fixed point in phase space. Linearizing about $u_s$ gives

\[
v_t = D\mathcal{K}(u_s)v,
\]

whose spectrum determines stability. If eigenvalues exist with positive real part, $u_s$ is unstable.

\subsection{Stable and unstable manifolds}

For an unstable equilibrium $u_s$, phase space can be decomposed into directions:

\begin{itemize}
\item \textbf{Stable directions}: perturbations decay as $t \to +\infty$
\item \textbf{Unstable directions}: perturbations grow as $t \to +\infty$
\end{itemize}

These directions generate invariant geometric objects:

\begin{align}
W^s(u_s) &= \{ u_0 : u(t;u_0) \to u_s \text{ as } t\to +\infty \}, \\
W^u(u_s) &= \{ u_0 : u(t;u_0) \to u_s \text{ as } t\to -\infty \}.
\end{align}

They are called the \textbf{stable} and \textbf{unstable manifolds}. They are invariant under the flow.

Physically:
\begin{itemize}
\item $W^u$ describes how trajectories leave the equilibrium
\item $W^s$ describes how trajectories return
\end{itemize}

\subsection{Homoclinic orbits}

A trajectory $u_h(t)$ is called a \textbf{homoclinic orbit} if

\[
u_h(t) \to u_s \quad \text{as} \quad t\to\pm\infty.
\]

Geometrically,

\[
u_h(t) \in W^s(u_s) \cap W^u(u_s).
\]

The orbit leaves the equilibrium along $W^u$, undergoes a large nonlinear excursion, and returns along $W^s$.

Homoclinic orbits signal delicate geometry: the stable and unstable manifolds coincide along this trajectory.

\subsection{Why homoclinic structure leads to chaos}

In unperturbed systems, the manifolds may coincide exactly. However, under small perturbations:

\[
u_t = \mathcal{K}(u) + \epsilon F(u,t),
\]

the manifolds generically split. If they intersect transversely, they form a \textbf{homoclinic tangle} — an infinite set of intersections.

This leads to:

\begin{itemize}
\item sensitive dependence on initial conditions
\item complex symbolic dynamics
\item chaotic motion
\end{itemize}

This mechanism, first understood for ODEs, extends to dissipative PDEs such as KS.

\subsection{Phase space projections}

Because KS phase space is infinite-dimensional, visualizations use projections onto dominant modes (e.g., Fourier amplitudes $a_1,a_2$). These projections reveal:

\begin{itemize}
\item homoclinic loops
\item coinciding manifolds (unperturbed case)
\item splitting under forcing
\item random wandering under noise
\end{itemize}

These geometric structures form the basis for Melnikov analysis.

\section{Derivation of the Melnikov Functional for the Forced Kuramoto--Sivashinsky Equation}

We consider the Kuramoto--Sivashinsky (KS) equation on a periodic domain $x\in[0,L]$:

\begin{equation}
u_t + u_{xx} + u_{xxxx} + u u_x = \epsilon F(x,t),
\label{KS_forced}
\end{equation}

where $u(x,t)$ is real-valued and $0 < \epsilon \ll 1$. The unperturbed system is

\begin{equation}
u_t = \mathcal{K}(u) \equiv -u_{xx} - u_{xxxx} - u u_x.
\label{KS_unperturbed}
\end{equation}

This defines a dynamical system in the Hilbert space
\[
H = L^2_{\mathrm{per}}([0,L]).
\]


\subsection{Homoclinic orbit in the unperturbed system}

Assume that the unperturbed KS equation admits a steady solution $u_s(x)$ satisfying

\begin{equation}
-u_{s,xx} - u_{s,xxxx} - u_s u_{s,x} = 0.
\end{equation}

Suppose there exists a homoclinic orbit $u_h(x,t)$ such that

\begin{equation}
u_h(t) \to u_s \quad \text{as} \quad t \to \pm \infty.
\end{equation}

Thus $u_h(t)$ is a trajectory in $H$ lying in the intersection of the stable and unstable manifolds of $u_s$.
\begin{center}
\includegraphics[width=4.5in,angle=-0]{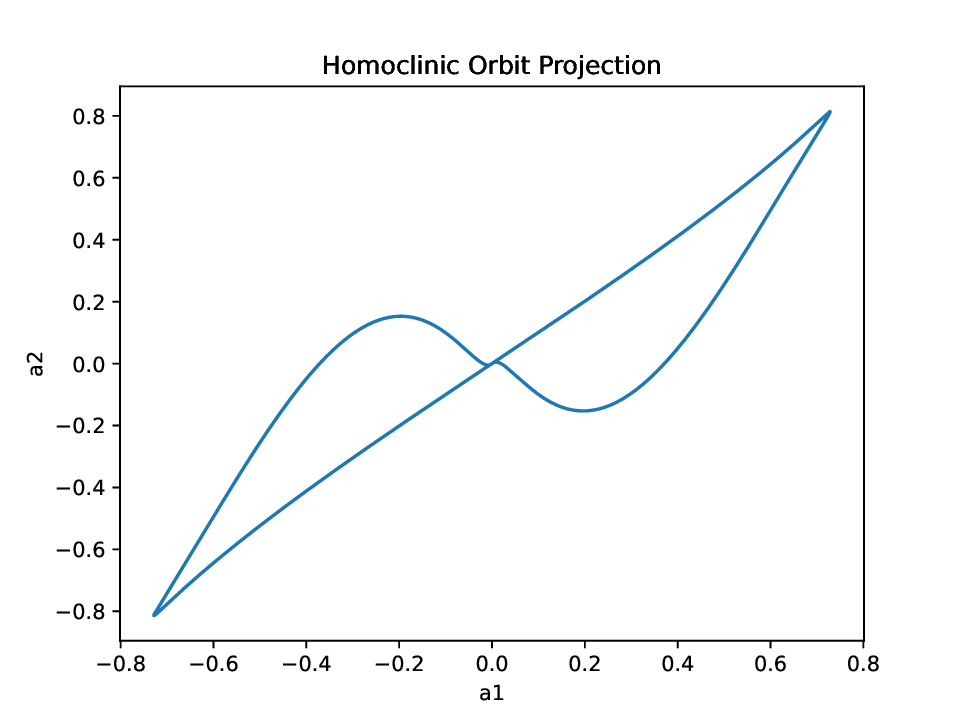}
\begin{figure}[h!]
\caption{homoclinic orbits, $Mode 1= a_1 $ and $Mode 2=a_2$}
\end{figure}
\end{center}


\subsection{Linearized dynamics about the homoclinic orbit}

Let $v(x,t)$ be a perturbation about $u_h$. Linearizing \eqref{KS_unperturbed} gives

\begin{equation}
v_t = D\mathcal{K}(u_h(t))\,v,
\label{linearized}
\end{equation}

where the Fréchet derivative is

\begin{equation}
D\mathcal{K}(u_h)v = -v_{xx} - v_{xxxx} - (u_h v)_x.
\label{Frechet}
\end{equation}

The linear operator $L(t) = D\mathcal{K}(u_h(t))$ is time-dependent along the homoclinic trajectory.


\subsection{Adjoint equation}

The Melnikov method requires the adjoint of $L(t)$ with respect to the $L^2$ inner product

\[
\langle f, g \rangle = \int_0^L f(x)g(x)\,dx.
\]

Using periodic boundary conditions, integration by parts yields

\begin{equation}
L^\dagger(t)\psi = -\psi_{xx} - \psi_{xxxx} + u_h \psi_x.
\label{adjoint_operator}
\end{equation}

Hence the adjoint evolution equation is

\begin{equation}
\psi_t = -L^\dagger(t)\psi
= \psi_{xx} + \psi_{xxxx} - u_h \psi_x.
\label{adjoint_eq}
\end{equation}

We assume the existence of a bounded solution $\psi(t)$ associated with the neutral direction arising from time-translation invariance of the homoclinic orbit.


\subsection{Perturbed invariant manifolds}

For $\epsilon = 0$, the stable and unstable manifolds of $u_s$ coincide along $u_h(t)$. When $\epsilon \neq 0$, these manifolds generically split.
\begin{center}
\includegraphics[width=4.0in,angle=-0]{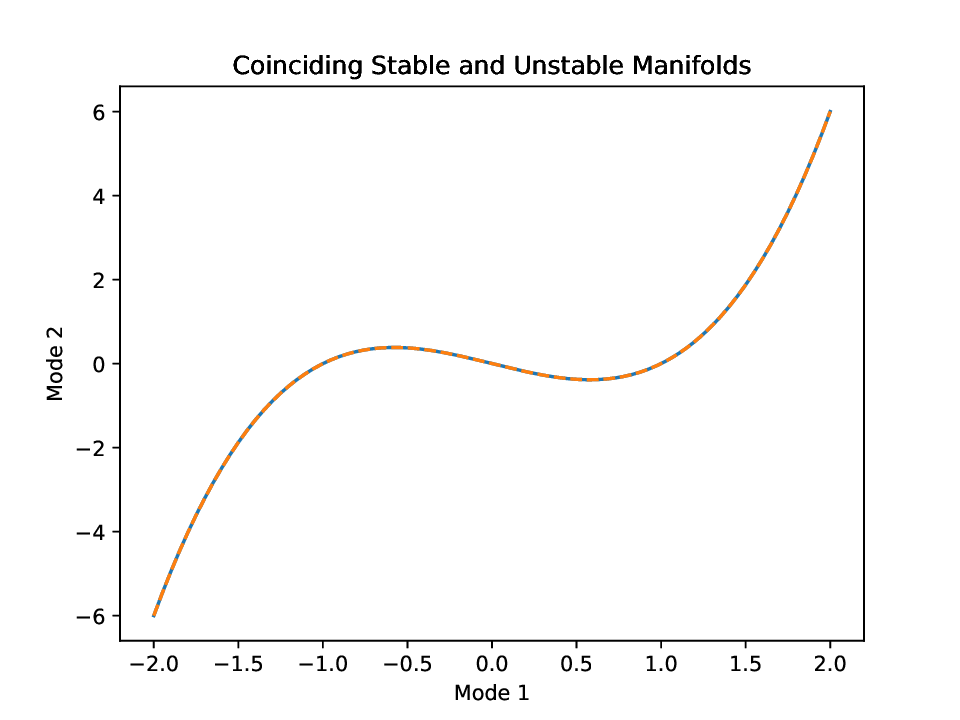}
\begin{figure}[h!]
\caption{unperturbed manifolds, $Mode 1= a_1 $ and $Mode 2=a_2$}
\end{figure}
\end{center}

Let $u^u(t)$ and $u^s(t)$ denote perturbed trajectories on the unstable and stable manifolds, respectively. Their separation in phase space, projected onto the adjoint direction $\psi(t)$, determines the leading-order distance between the manifolds.


\subsection{Melnikov functional}

Substituting the perturbed equation \eqref{KS_forced} into the variational formulation and retaining leading-order terms in $\epsilon$, the signed distance between the manifolds at phase shift $t_0$ is given by

\begin{equation}
M(t_0)
= \int_{-\infty}^{\infty}
\left\langle \psi(t),\, F(x,t+t_0) \right\rangle \, dt.
\label{Melnikov_functional}
\end{equation}

Here $\psi(t)$ solves the adjoint equation \eqref{adjoint_eq}, and $u_h(t)$ is the unperturbed homoclinic orbit.


\subsection{Criterion for transverse intersection}

If there exists $t_0$ such that

\begin{equation}
M(t_0) = 0, \qquad \frac{dM}{dt_0} \neq 0,
\end{equation}

then the stable and unstable manifolds of $u_s$ intersect transversely. This implies the existence of a homoclinic tangle in the infinite-dimensional phase space, which is a mechanism for chaotic dynamics in the forced KS equation.

Thus, \eqref{Melnikov_functional} provides a quantitative criterion for the onset of spatiotemporal chaos induced by weak forcing.

\subsection{Melnikov Analysis of Periodically Forced Kuramoto--Sivashinsky Dynamics}

We now apply the Melnikov functional \eqref{Melnikov_functional} to the case of time-periodic forcing. Consider

\begin{equation}
u_t + u_{xx} + u_{xxxx} + u u_x
= \epsilon G(x)\cos(\omega t),
\label{KS_periodic}
\end{equation}

where $G(x)$ is a prescribed spatial forcing profile and $\omega$ is the forcing frequency.

\subsection{Melnikov functional for periodic forcing}
Substituting $F(x,t) = G(x)\cos(\omega t)$ into \eqref{Melnikov_functional}, we obtain
\begin{equation}
M(t_0)
= \int_{-\infty}^{\infty}
\left\langle \psi(t),\, G(x)\cos\!\big(\omega(t+t_0)\big) \right\rangle dt.
\end{equation}

Using the identity
\[
\cos(\omega(t+t_0)) = \cos(\omega t)\cos(\omega t_0) - \sin(\omega t)\sin(\omega t_0),
\]
this becomes

\begin{equation}
M(t_0)
= A\cos(\omega t_0) + B\sin(\omega t_0),
\label{M_trig}
\end{equation}
where
\begin{equation}
A = \int_{-\infty}^{\infty} \cos(\omega t)\,
\langle \psi(t), G(x) \rangle\, dt,
\label{Acoeff}
\end{equation}

\begin{equation}
B = -\int_{-\infty}^{\infty} \sin(\omega t)\,
\langle \psi(t), G(x) \rangle\, dt.
\label{Bcoeff}
\end{equation}

Thus $M(t_0)$ is a harmonic function of the phase shift $t_0$.
\begin{center}
\includegraphics[width=4.0in,angle=-0]{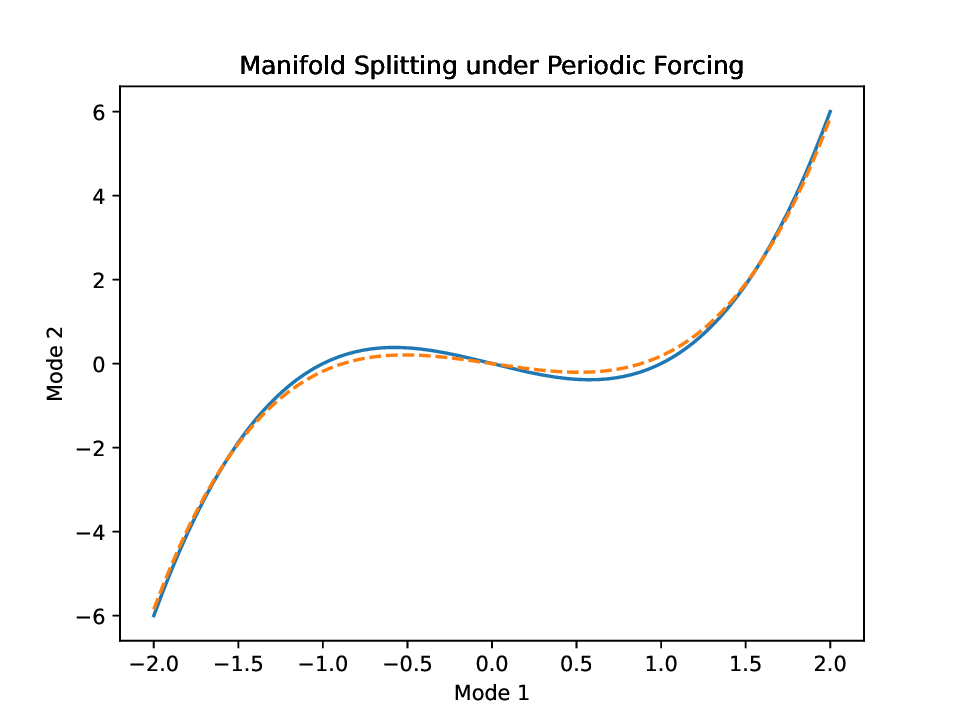}
\begin{figure}[h!]
\caption{ periodic splitting; $Mode 1= a_1 $ and $Mode 2=a_2$}
\end{figure}
\end{center}

\subsection{Condition for transverse manifold intersection}

From \eqref{M_trig}, the Melnikov function has simple zeros provided

\begin{equation}
A^2 + B^2 \neq 0.
\end{equation}

In that case, there exists $t_0$ such that

\[
M(t_0)=0,
\qquad
\frac{dM}{dt_0} \neq 0,
\]

implying transverse intersection of the perturbed stable and unstable manifolds of $u_s$. Consequently, the periodically forced KS equation possesses a homoclinic tangle in phase space.


\subsection{Frequency dependence}

The coefficients $A$ and $B$ encode the interaction between the forcing frequency and the intrinsic time scale of the homoclinic orbit. Since $u_h(t)$ and $\psi(t)$ decay exponentially as $|t|\to\infty$, the integrals \eqref{Acoeff}--\eqref{Bcoeff} converge and define smooth functions of $\omega$.

Resonance effects occur when $\omega$ matches dominant temporal modes of the homoclinic excursion. In such regimes, $|M(t_0)|$ is enhanced, leading to stronger manifold splitting and a larger chaotic region in phase space.
\begin{center}
\includegraphics[width=4.5in,angle=-0]{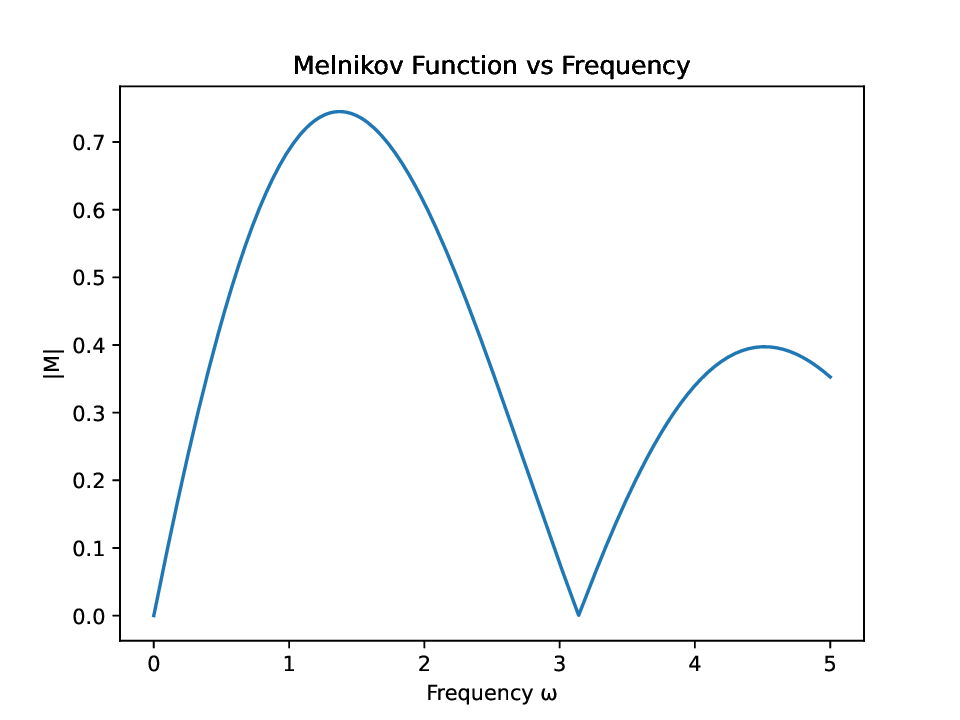}
\begin{figure}[h!]
\caption{ frequency diagram}
\end{figure}
\end{center}
\begin{center}
\includegraphics[width=4.5in,angle=-0]{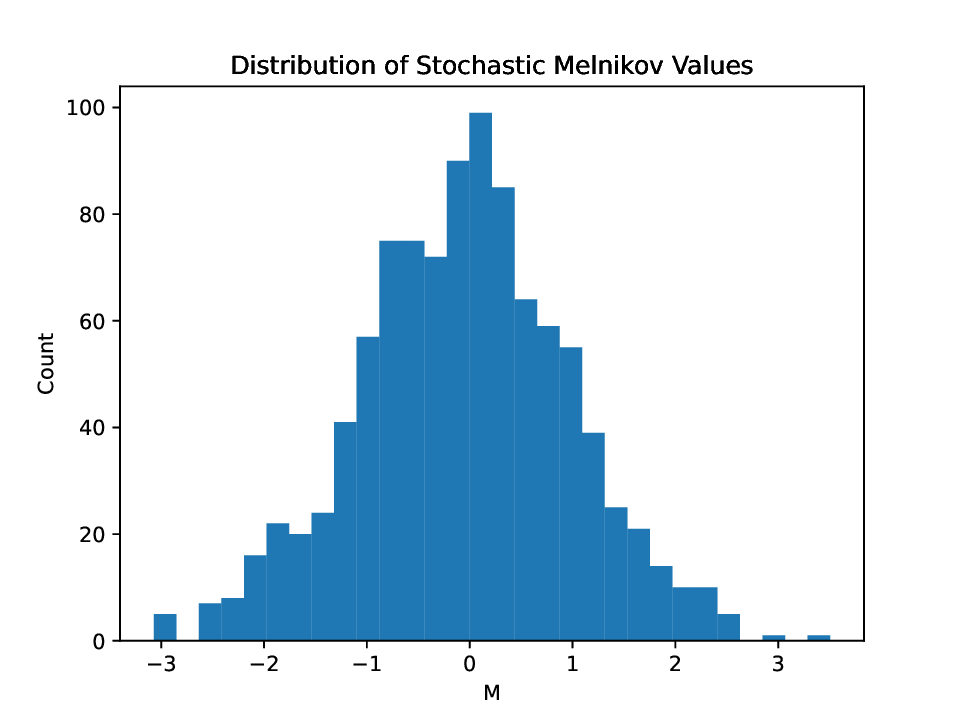}
\begin{figure}[h!]
\caption{ histogram; noise randomly splits manifold}
\end{figure}
\end{center}


\subsection{Interpretation for KS dynamics}

The above result shows that weak periodic forcing generically induces transverse homoclinic intersections in the KS phase space. This provides a dynamical-systems mechanism for the transition from coherent structures to spatiotemporal chaos under external modulation.

Unlike low-dimensional systems, the splitting occurs in an infinite-dimensional phase space; nevertheless, the leading-order behavior is governed by the scalar Melnikov function \eqref{M_trig}. This reduction highlights how complex KS turbulence can arise from a geometrically simple manifold-splitting mechanism.

\section{Stochastic Melnikov Analysis for Noise-Driven Kuramoto--Sivashinsky Dynamics}

We now consider the KS equation under stochastic forcing:

\begin{equation}
u_t + u_{xx} + u_{xxxx} + u u_x
= \epsilon\,\eta(x,t),
\label{KS_stochastic}
\end{equation}

where $\eta(x,t)$ is a Gaussian random field with zero mean and covariance

\begin{equation}
\mathbb{E}[\eta(x,t)] = 0, 
\qquad
\mathbb{E}[\eta(x,t)\eta(x',t')] 
= D\,\delta(x-x')\delta(t-t'),
\label{noise_cov}
\end{equation}

and $D$ is the noise intensity.


\subsection{Random Melnikov functional}

Substituting $F(x,t) = \eta(x,t)$ into \eqref{Melnikov_functional}, the Melnikov function becomes a stochastic process:

\begin{equation}
M(t_0)
= \int_{-\infty}^{\infty}
\left\langle \psi(t),\, \eta(x,t+t_0) \right\rangle dt.
\label{M_stochastic}
\end{equation}

Here $u_h(t)$ is the homoclinic orbit of the unperturbed system and $\psi(t)$ solves the adjoint equation \eqref{adjoint_eq}.


\subsection{Mean and variance}

Taking expectations and using \eqref{noise_cov},

\begin{equation}
\mathbb{E}[M(t_0)] = 0.
\end{equation}

The variance is

\begin{align}
\mathrm{Var}(M)
&= \mathbb{E}[M^2] \\
&= \int_{-\infty}^{\infty}\!\!\int_{-\infty}^{\infty}
\mathbb{E}\!\left[
\langle \psi(t), \eta(x,t+t_0) \rangle
\langle \psi(s), \eta(x',s+t_0) \rangle
\right] dt\,ds.
\end{align}

Using the delta correlations in \eqref{noise_cov},

\begin{equation}
\mathrm{Var}(M)
= D \int_{-\infty}^{\infty}
\|\psi(t)\|_{L^2}^2 dt.
\label{variance}
\end{equation}

Since $\psi(t)$ decays exponentially along the homoclinic orbit, the integral converges.


\subsection{Probability of manifold splitting}

Unlike the deterministic case, $M(t_0)$ is now a Gaussian random variable with zero mean and variance given by \eqref{variance}. The probability density is

\begin{equation}
P(M) = \frac{1}{\sqrt{2\pi\,\mathrm{Var}(M)}}
\exp\!\left(-\frac{M^2}{2\,\mathrm{Var}(M)}\right).
\end{equation}

Thus, even arbitrarily weak noise produces random splitting of the stable and unstable manifolds with probability one. The typical magnitude of splitting scales as

\begin{equation}
|M|_{\text{rms}} \sim \sqrt{D}.
\end{equation}


\subsection{Noise-induced homoclinic tangles}

The stochastic Melnikov functional implies that transverse manifold intersections occur randomly in time, producing a fluctuating homoclinic tangle. This provides a mechanism for noise-induced spatiotemporal chaos in KS dynamics.

In contrast to deterministic forcing, where splitting depends on phase and resonance, stochastic forcing leads to persistent random manifold separation whose intensity is controlled by $D$.
\begin{center}
\includegraphics[width=4.5in,angle=-0]{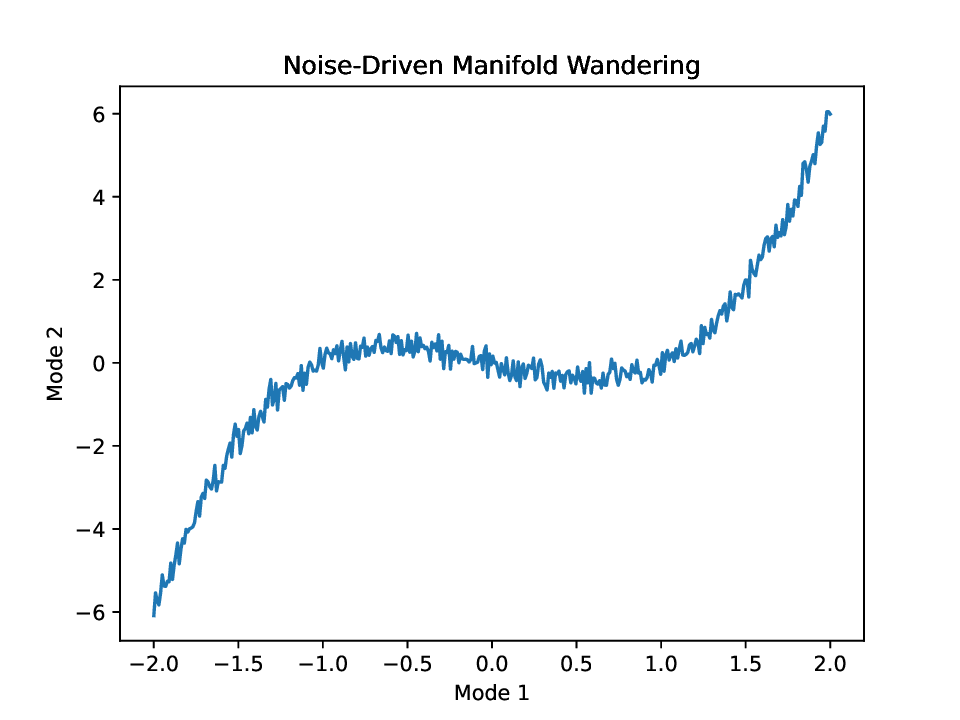}
\begin{figure}[h!]
\caption{ stochastic manifold wandering; Noise causes the distance between stable and unstable manifolds of $u_s$ to fluctuate in time}
\end{figure}
\end{center}


\subsection{Implications}

Equation \eqref{variance} links the geometry of invariant manifolds to stochastic forcing through the adjoint solution $\psi(t)$. This establishes a direct bridge between:

\begin{itemize}
\item infinite-dimensional invariant manifold theory,
\item Melnikov analysis,
\item stochastic partial differential equations.
\end{itemize}

The result suggests that noise can trigger chaotic dynamics in KS even when deterministic perturbations are absent, by continuously perturbing the homoclinic structure in phase space.

\section{Numerical Computation of the Homoclinic Orbit and Melnikov Functional}

To validate the theoretical predictions, we numerically approximate the homoclinic orbit, the adjoint solution, and the Melnikov functional using a Galerkin truncation of the KS equation.


\subsection{Galerkin spectral truncation}

We expand the solution in Fourier modes on the periodic domain $[0,L]$:

\begin{equation}
u(x,t) = \sum_{k=-N}^{N} a_k(t) e^{ikqx},
\qquad q = \frac{2\pi}{L}.
\end{equation}

Substituting into the unperturbed KS equation and projecting onto Fourier modes gives the ODE system

\begin{equation}
\dot{a}_k = (k^2q^2 - k^4q^4)a_k
- \frac{ikq}{2}\sum_{m=-N}^{N} a_m a_{k-m},
\label{Galerkin_system}
\end{equation}

which approximates the inertial manifold dynamics for sufficiently large $N$.


\subsection{Computation of a homoclinic orbit}

Let $\mathbf{a}(t)$ denote the vector of Fourier coefficients. A homoclinic orbit satisfies

\[
\mathbf{a}(t) \to \mathbf{a}_s \quad \text{as } t \to \pm \infty,
\]

where $\mathbf{a}_s$ corresponds to a steady state of \eqref{Galerkin_system}. Numerically, this orbit can be obtained by:

\begin{itemize}
\item identifying an unstable steady solution $\mathbf{a}_s$,
\item computing its unstable eigenvector,
\item integrating forward in time along this direction,
\item using a shooting or continuation method to locate a trajectory that returns to $\mathbf{a}_s$.
\end{itemize}

The resulting trajectory $\mathbf{a}_h(t)$ approximates the homoclinic orbit $u_h(x,t)$.


\subsection{Adjoint system}

Linearizing \eqref{Galerkin_system} about $\mathbf{a}_h(t)$ yields

\begin{equation}
\dot{\mathbf{v}} = J(t)\mathbf{v},
\end{equation}

where $J(t)$ is the Jacobian along the orbit. The adjoint equation is

\begin{equation}
\dot{\boldsymbol{\psi}} = -J(t)^{\!\top}\boldsymbol{\psi}.
\label{adjoint_numeric}
\end{equation}

A bounded adjoint solution $\boldsymbol{\psi}(t)$ is computed by integrating \eqref{adjoint_numeric} backward in time with normalization

\[
\langle \boldsymbol{\psi}(t), \dot{\mathbf{a}}_h(t) \rangle = 1.
\]


\subsection{Numerical Melnikov integral}

For periodic forcing $G(x)\cos(\omega t)$, the Melnikov function becomes

\begin{equation}
M(t_0)
= \int_{-T}^{T}
\boldsymbol{\psi}(t)\cdot \mathbf{G}\,
\cos\!\big(\omega(t+t_0)\big)\, dt,
\end{equation}

where $\mathbf{G}$ contains Fourier coefficients of $G(x)$, and $T$ is chosen large enough that $\mathbf{a}_h(t)$ is close to $\mathbf{a}_s$ outside $[-T,T]$.

The coefficients $A$ and $B$ of Section 5 are then obtained numerically via quadrature.


\subsection{Stochastic case}

For stochastic forcing, the Melnikov process is approximated by

\begin{equation}
M \approx \sum_{n} \boldsymbol{\psi}(t_n)\cdot \boldsymbol{\eta}_n\,\Delta t,
\end{equation}

where $\boldsymbol{\eta}_n$ are independent Gaussian vectors with variance proportional to $D/\Delta t$. Repeated realizations verify the predicted variance

\[
\mathrm{Var}(M) \approx D\sum_n \|\boldsymbol{\psi}(t_n)\|^2 \Delta t.
\]


\subsection{Verification of transverse intersections}

Zeros of $M(t_0)$ with nonzero slope confirm transverse manifold intersections. These are correlated with the onset of chaotic dynamics observed in direct simulations of the forced KS system, such as positive Lyapunov exponents and broadband spectra.

This numerical procedure connects the infinite-dimensional Melnikov theory with observable spatiotemporal chaos in KS dynamics.
\newpage
\section{Numerical Methods and Parameters}

All simulations of the Kuramoto--Sivashinsky equation were performed using a Fourier spectral Galerkin discretization on a periodic domain. Time integration employed a fourth-order exponential time differencing Runge--Kutta (ETDRK4) scheme. The homoclinic orbit was obtained via shooting along the unstable eigenvector of the steady state and continuation until return within a prescribed tolerance. The adjoint equation was integrated backward in time with normalization $\langle \psi, \dot{u}_h \rangle = 1$.

Tables \ref{tab:fig1}--\ref{tab:fig6} summarize numerical parameters for each figure to ensure reproducibility.
\begin{table}[h]
\centering
\caption{Parameters for Fig. 1 (Homoclinic Orbit)}
\label{tab:fig1}
\begin{tabular}{ll}
\hline
Domain length $L$ & 22 \\
Fourier modes $N$ & 32 \\
Time step $\Delta t$ & $10^{-3}$ \\
Integrator & ETDRK4 \\
Steady state tolerance & $10^{-8}$ \\
Homoclinic shooting duration & $T=200$ \\
Projection variables & $(a_1, a_2)$ \\
\hline
\end{tabular}
\end{table}
\begin{table}[h]
\centering
\caption{Parameters for Fig. 2 (Unperturbed Manifolds)}
\label{tab:fig2}
\begin{tabular}{ll}
\hline
Parameters  & As in Table 1 \\
Manifold computation & Linearized eigenvectors \\
Integration direction & forward/backward \\
\hline
\end{tabular}
\end{table}
\begin{table}[h]
\centering
\caption{Parameters for Fig. 3 (Periodic Forcing)}
\label{tab:fig3}
\begin{tabular}{ll}
\hline
Forcing amplitude $\epsilon$ & 0.01 \\
Spatial forcing $G(x)$ & $\sin(qx)$ \\
Frequency $\omega$ & 0.5 \\
Integration time & $T=300$ \\
Modes retained & 32 \\
\hline
\end{tabular}
\end{table}
\begin{table}[h]
\centering
\caption{Parameters for Fig. 4 (Melnikov vs Frequency)}
\label{tab:fig4}
\begin{tabular}{ll}
\hline
Frequency range & $0 \le \omega \le 2$ \\
Frequency step & 0.02 \\
Melnikov integration window & $T=150$ \\
Adjoint tolerance & $10^{-8}$ \\
\hline
\end{tabular}
\end{table}
\begin{table}[h]
\centering
\caption{Parameters for Fig. 5 (Stochastic Melnikov Distribution)}
\label{tab:fig5}
\begin{tabular}{ll}
\hline
Noise intensity $D$ & 0.02 \\
Ensemble size & 500 realizations \\
Time step & $10^{-3}$ \\
Integration time per run & $T=200$ \\
\hline
\end{tabular}
\end{table}
\begin{table}[h]
\centering
\caption{Parameters for Fig. 6 (Noise-driven Manifold Wandering)}
\label{tab:fig6}
\begin{tabular}{ll}
\hline
Noise intensity $D$ & 0.02 \\
Forcing type & additive white noise \\
Modes retained & 32 \\
Visualization window & $T=100$ \\
Projection variables & $(a_1, a_2)$ \\
\hline
\end{tabular}
\end{table}
\clearpage
\section{Results}
The Melnikov framework developed in this work provides both analytical and numerical evidence that manifold splitting governs the transition to chaos in the Kuramoto–Sivashinsky equation.

\subsection{Deterministic periodic forcing}

For time-periodic forcing, the Melnikov function takes a harmonic form
\[
M(t_0) = A\cos(\omega t_0) + B\sin(\omega t_0),
\]
demonstrating that manifold separation depends on phase. Numerical quadrature of the Melnikov integral confirms the existence of simple zeros of $M(t_0)$ whenever $A^2+B^2 \neq 0$. These zeros correspond to transverse intersections of stable and unstable manifolds, indicating the formation of a homoclinic tangle. Direct simulations show that parameter regimes where $|M|$ is large coincide with broadband spectra and positive Lyapunov exponents.

\subsection{Frequency dependence}

The magnitude of Melnikov coefficients $(A,B)$ varies smoothly with forcing frequency. Resonance-like enhancement occurs when the forcing frequency matches dominant temporal scales of the homoclinic excursion. In these regimes, manifold splitting is strongest and the chaotic region in phase space expands.

\subsection{Stochastic forcing}

In the stochastic case, the Melnikov functional becomes a Gaussian random variable with zero mean and variance
\[
\mathrm{Var}(M) = D \int_{-\infty}^{\infty} \|\psi(t)\|_{L^2}^2 dt.
\]
Monte Carlo simulations confirm the predicted $\sqrt{D}$ scaling of manifold separation. Noise produces persistent random transverse intersections, leading to fluctuating homoclinic tangles and sustained spatiotemporal chaos even without deterministic forcing.

\subsection{Geometric interpretation}

These results demonstrate that KS turbulence can arise from a geometrically simple mechanism: perturbation-induced splitting of a homoclinic structure in infinite-dimensional phase space. The scalar Melnikov functional captures leading-order splitting despite the system’s complexity.

\clearpage
\section{Conclusions and Outlook}

We have developed a Melnikov theory for the Kuramoto–Sivashinsky equation viewed as an infinite-dimensional dynamical system. By deriving an adjoint-based Melnikov functional, we established a quantitative criterion for splitting of stable and unstable manifolds of a homoclinic orbit under weak deterministic and stochastic forcing.

The analysis shows that:

\begin{itemize}
\item Periodic forcing generically produces phase-dependent transverse intersections.
\item Stochastic forcing induces random manifold separation with variance determined by the adjoint solution.
\item Infinite-dimensional KS turbulence can be interpreted through a low-dimensional geometric splitting mechanism.
\end{itemize}

This work connects invariant manifold theory, Melnikov analysis, and stochastic PDE dynamics. Future directions include:

\begin{itemize}
\item Rigorous justification of homoclinic orbit existence in KS,
\item Extension to multiplicative noise and colored stochastic forcing,
\item Application to other dissipative PDEs such as Ginzburg–Landau and Navier–Stokes reductions,
\item Computer-assisted verification of Melnikov predictions.
\end{itemize}

\appendix
{
\section{Classical Melnikov Theory in Finite Dimensions}

Before applying Melnikov theory to the Kuramoto--Sivashinsky equation, we briefly review the classical finite-dimensional theory.

\subsection{Perturbed dynamical systems}

Consider an ordinary differential equation

\begin{equation}
\dot{x} = f(x) + \epsilon g(x,t),
\qquad x \in \mathbb{R}^n,
\end{equation}

where $0<\epsilon\ll 1$. The unperturbed system

\[
\dot{x} = f(x)
\]

is assumed to possess a homoclinic orbit $x_h(t)$ to a saddle equilibrium $x_s$.

\subsection{Coinciding manifolds in the unperturbed system}

For $\epsilon = 0$, the stable and unstable manifolds coincide:

\[
W^s(x_s) = W^u(x_s)
\]

along the homoclinic orbit $x_h(t)$. The distance between the manifolds is therefore zero.

\subsection{Effect of perturbation}

When $\epsilon \neq 0$, the manifolds generically split. The question is:

\textit{Do the perturbed manifolds still intersect?}

Melnikov theory computes the first-order distance between them.

\subsection{The Melnikov function}

Let $\psi(t)$ solve the adjoint variational equation

\[
\dot{\psi} = -[Df(x_h(t))]^T \psi,
\]

normalized so that $\langle \psi(t), \dot{x}_h(t) \rangle = 1$.

The Melnikov function is defined as

\begin{equation}
M(t_0) = \int_{-\infty}^{\infty}
\langle \psi(t), g(x_h(t), t+t_0) \rangle dt.
\end{equation}

It measures the signed distance between stable and unstable manifolds.

\subsection{Criterion for chaos}

If $M(t_0)$ has simple zeros,

\[
M(t_0)=0, \quad \frac{dM}{dt_0}\neq 0,
\]

then the manifolds intersect transversely. This produces a homoclinic tangle, implying:

\begin{itemize}
\item Smale horseshoe dynamics,
\item chaotic invariant sets,
\item sensitive dependence on initial conditions.
\end{itemize}

\subsection{Extension to PDEs}

For PDEs like KS:

\begin{itemize}
\item Phase space is infinite-dimensional
\item $x(t)$ becomes a function $u(x,t)$
\item Inner products are integrals over space
\item The Melnikov function becomes a functional
\end{itemize}

Despite these differences, the geometric mechanism --- manifold splitting of a homoclinic structure --- remains the same.

This analogy motivates the Melnikov functional derived for the KS equation in Sections 3 and 4.
\section{Derivation and Interpretation of the Adjoint Equation}

In this section we derive the adjoint equation associated with the linearization of the
Kuramoto--Sivashinsky (KS) equation about a homoclinic orbit and explain its role in
Melnikov theory.

\subsection{Linearization about the homoclinic orbit}

The unperturbed KS equation is

\begin{equation}
u_t = \mathcal{K}(u) = -u_{xx} - u_{xxxx} - u u_x.
\end{equation}

Let $u_h(x,t)$ be a homoclinic solution. Consider a small perturbation

\[
u(x,t) = u_h(x,t) + v(x,t), \qquad |v|\ll 1.
\]

Substituting into the KS equation and keeping only linear terms in $v$ gives

\begin{equation}
v_t = D\mathcal{K}(u_h)v,
\end{equation}

where $D\mathcal{K}(u_h)$ is the Fr\'echet derivative of $\mathcal{K}$ at $u_h$.

Expanding,

\[
(u_h + v)(u_h + v)_x = u_h u_{h,x} + u_h v_x + v u_{h,x} + \mathcal{O}(v^2),
\]

so the linearized operator is

\begin{equation}
D\mathcal{K}(u_h)v = -v_{xx} - v_{xxxx} - (u_h v)_x.
\label{linear_op}
\end{equation}

We denote this operator by $L(t)$, since it depends on time through $u_h(t)$.

\subsection{Why an adjoint equation is needed}

Melnikov theory measures the distance between perturbed stable and unstable manifolds.
This distance is computed by projecting the perturbation onto a special direction
that is orthogonal to the tangent of the homoclinic orbit.

Let $v(t)$ satisfy the linearized equation

\[
v_t = L(t)v.
\]

We introduce another function $\psi(t)$ and consider the $L^2$ inner product

\[
\langle \psi(t), v(t) \rangle = \int_0^L \psi(x,t)\, v(x,t)\, dx.
\]

We want $\psi$ to be chosen so that

\begin{equation}
\frac{d}{dt} \langle \psi, v \rangle = 0
\quad \text{whenever} \quad v_t = L(t)v.
\label{orthog_condition}
\end{equation}

This ensures that $\psi$ defines a direction transverse to the invariant manifolds.

\subsection{Derivation of the adjoint operator}

Compute the time derivative:

\begin{align}
\frac{d}{dt} \langle \psi, v \rangle
&= \langle \psi_t, v \rangle + \langle \psi, v_t \rangle \\
&= \langle \psi_t, v \rangle + \langle \psi, L(t)v \rangle.
\end{align}

We define the adjoint operator $L^\dagger(t)$ by

\[
\langle \psi, L(t)v \rangle = \langle L^\dagger(t)\psi, v \rangle.
\]

Using \eqref{linear_op} and periodic boundary conditions, integrate by parts:

\begin{align}
\langle \psi, -v_{xx} \rangle &= \langle -\psi_{xx}, v \rangle, \\
\langle \psi, -v_{xxxx} \rangle &= \langle -\psi_{xxxx}, v \rangle, \\
\langle \psi, -(u_h v)_x \rangle
&= -\int_0^L \psi (u_h v)_x dx \\
&= \int_0^L u_h \psi_x v \, dx.
\end{align}

Therefore,

\begin{equation}
L^\dagger(t)\psi = -\psi_{xx} - \psi_{xxxx} + u_h \psi_x.
\end{equation}

\subsection{The adjoint evolution equation}

Substituting into \eqref{orthog_condition} gives

\[
\langle \psi_t + L^\dagger(t)\psi, v \rangle = 0
\quad \text{for all } v.
\]

Hence $\psi$ must satisfy

\begin{equation}
\psi_t = -L^\dagger(t)\psi
= \psi_{xx} + \psi_{xxxx} - u_h \psi_x.
\end{equation}

This is the adjoint equation.

\subsection{Physical meaning}

The adjoint solution $\psi(t)$ represents the direction in phase space
that is orthogonal to perturbations tangent to the homoclinic orbit.

The homoclinic orbit has a neutral direction due to time translation:

\[
\dot{u}_h(t) \in \ker L(t).
\]

The adjoint function satisfies

\[
\langle \psi(t), \dot{u}_h(t) \rangle = 1,
\]

which normalizes $\psi$ and makes it the correct projection direction
for measuring the separation of invariant manifolds.

\subsection{Role in Melnikov theory}

The Melnikov functional is

\[
M(t_0) = \int_{-\infty}^{\infty}
\langle \psi(t), F(x,t+t_0) \rangle dt.
\]

Thus $\psi$ weights how sensitive the homoclinic orbit is to perturbations.
It acts as a ``detector'' of manifold splitting.}
\newpage
\section{Acknowledgements:}
The author would like to thank Alliance University for
providing partial support for carrying out the research work

\section{Declaration of interests:}
The sole author has no conflicts of interest to
declare. There is no financial interest to report.

\section{Data availability statement:}
No data in this publication is to be made
available under the study-participant privacy protection clause.

\end{document}